\documentclass[11pt]{amsart}   

\usepackage[latin1]{inputenc}
\usepackage[T1]{fontenc}
\usepackage[english]{babel}

\usepackage{amsmath}
\usepackage{amssymb}
\usepackage{amsthm}
\usepackage{mathabx}
\usepackage{MnSymbol}

\usepackage{color}
\usepackage{tikz}
\usepackage{multicol}

\theoremstyle{plain}
\newtheorem{theorem}{Theorem}
\newtheorem{lemma}[theorem]{Lemma}

\newtheorem{corollary}[theorem]{Corollary}

\newtheorem{example}{\it Example\/}
\author[L. Colmenarejo and M. Rosas]{Laura Colmenarejo and Mercedes Rosas} 
\thanks{Both authors are partially supported by MTM2013-40455-P, P12-FQM-2696, FQM-333, and FEDER. }

\title[Combinatorics on reduced Kronecker coefficients]{\textsc{Combinatorics on a family of \\
reduced Kronecker coefficients}}
\email{laurach@us.es, mrosas@us.es}
\address{Department of Algebra, University of Seville}
\begin{document}

\begin{abstract}
{\it The reduced Kronecker coefficients are particular instances of Kronecker coefficients that contain enough information to recover them. In this notes we compute the generating function of a family of reduced Kronecker coefficients. We also gives its connection to the plane partitions, which allows us to check that this family satisfies the saturation conjecture for reduced Kronecker coefficients, and that they are weakly increasing. Thanks to its generating function we can describe our family by a quasipolynomial, specifying its degree and period. }

\end{abstract}

\maketitle
\vspace{0.7cm}

\section*{Introduction}

 With the original aim of trying to  understand the rate of grow experienced by the Kronecker coefficient as we increase the sizes of its rows, we investigate a  family of reduced Kronecker coefficients, $\overline{g}_{(k^a),(k^b)}^{(k)}$.

 The Kronecker coefficient are the multiplicities appearing in the decomposition into irreducible of the tensor product of two irreducible representations of the symmetric group.  They also appear naturally in the study of the general lineal group and the unitary group.
 The reduced Kronecker coefficients, on the other hand, are particular instances of Kronecker coefficients, believed to be simpler to understand, but that contain enough information to recover them. 

In Theorem \ref{ThmRowsKT} we compute the generating series  $\mathcal{F}_{a,b}=\sum_{k }\overline{g}_{(k^a),(k^b)}^{(k)}x^k$ for this family of coefficients. We achieve this by giving an explicit bijection between a family of coloured partitions, and the Kronecker tableaux of   Orellana and   Ballantine,  \cite{MR2264933}. Since $\mathcal{F}_{a,b}$  turns out to be a rational generating function, these coefficients obey a lineal recurrence. 

In Theorem \ref{ThmPPRows} we give an striking connection to plane partitions: The reduced Kronecker coefficient $\overline{g}_{(k^a),(k^a)}^{(k)}$ counts the number of plane partitions of $k$ fitting on a $2\times a$ rectangle.  We obtain this result by comparing the generating function obtained in Theorem \ref{ThmRowsKT}  with MacMahon's classical formula. 
\newpage 
For instance, the reduced Kronecker coefficient for $a=4$ and $k=3$ is $\overline{g}_{(3,3,3,3),(3,3,3,3)}^{(3)}=5$, and there are $5$ plane partitions of $3$ fitting inside an $2 \times 4$ rectangle\\
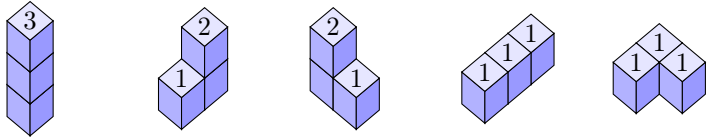
 \begin{center}
 {\small
  \begin{minipage}{0.15\linewidth}  
    \begin{tikzpicture}[x=(220:0.4cm), y=(-40:0.4cm), z=(90:0.4242cm)]
\foreach \m [count=\y] in {{3}}{
  \foreach \n [count=\x] in \m {
  \ifnum \n>0
      \foreach \z in {1,...,\n}{
        \draw [fill=blue!30] (\x+1,\y,\z) -- (\x+1,\y+1,\z) -- (\x+1, \y+1, \z-1) -- (\x+1, \y, \z-1) -- cycle;
        \draw [fill=blue!40] (\x,\y+1,\z) -- (\x+1,\y+1,\z) -- (\x+1, \y+1, \z-1) -- (\x, \y+1, \z-1) -- cycle;
        \draw [fill=blue!10] (\x,\y,\z)   -- (\x+1,\y,\z)   -- (\x+1, \y+1, \z)   -- (\x, \y+1, \z) -- cycle;  
      }
      \node at (\x+0.5, \y+0.5, \n) {\n};
 \fi
 }
}
\end{tikzpicture}
\end{minipage}
  \begin{minipage}{0.15\linewidth}  
    \begin{tikzpicture}[x=(220:0.4cm), y=(-40:0.4cm), z=(90:0.4242cm)]
\foreach \m [count=\y] in {{2,1}}{
  \foreach \n [count=\x] in \m {
  \ifnum \n>0
      \foreach \z in {1,...,\n}{
        \draw [fill=blue!30] (\x+1,\y,\z) -- (\x+1,\y+1,\z) -- (\x+1, \y+1, \z-1) -- (\x+1, \y, \z-1) -- cycle;
        \draw [fill=blue!40] (\x,\y+1,\z) -- (\x+1,\y+1,\z) -- (\x+1, \y+1, \z-1) -- (\x, \y+1, \z-1) -- cycle;
        \draw [fill=blue!10] (\x,\y,\z)   -- (\x+1,\y,\z)   -- (\x+1, \y+1, \z)   -- (\x, \y+1, \z) -- cycle;  
      }
      \node at (\x+0.5, \y+0.5, \n) {\n};
 \fi
 }
}
\end{tikzpicture}
\end{minipage}
  \begin{minipage}{0.15\linewidth}  
    \begin{tikzpicture}[x=(220:0.4cm), y=(-40:0.4cm), z=(90:0.4242cm)]
\foreach \m [count=\y] in {{2},{1}}{
  \foreach \n [count=\x] in \m {
  \ifnum \n>0
      \foreach \z in {1,...,\n}{
        \draw [fill=blue!30] (\x+1,\y,\z) -- (\x+1,\y+1,\z) -- (\x+1, \y+1, \z-1) -- (\x+1, \y, \z-1) -- cycle;
        \draw [fill=blue!40] (\x,\y+1,\z) -- (\x+1,\y+1,\z) -- (\x+1, \y+1, \z-1) -- (\x, \y+1, \z-1) -- cycle;
        \draw [fill=blue!10] (\x,\y,\z)   -- (\x+1,\y,\z)   -- (\x+1, \y+1, \z)   -- (\x, \y+1, \z) -- cycle;  
      }
      \node at (\x+0.5, \y+0.5, \n) {\n};
 \fi
 }
}
\end{tikzpicture}
\end{minipage}
  \begin{minipage}{0.15\linewidth}  
    \begin{tikzpicture}[x=(220:0.4cm), y=(-40:0.4cm), z=(90:0.4242cm)]
\foreach \m [count=\y] in {{1,1,1}}{
  \foreach \n [count=\x] in \m {
  \ifnum \n>0
      \foreach \z in {1,...,\n}{
        \draw [fill=blue!30] (\x+1,\y,\z) -- (\x+1,\y+1,\z) -- (\x+1, \y+1, \z-1) -- (\x+1, \y, \z-1) -- cycle;
        \draw [fill=blue!40] (\x,\y+1,\z) -- (\x+1,\y+1,\z) -- (\x+1, \y+1, \z-1) -- (\x, \y+1, \z-1) -- cycle;
        \draw [fill=blue!10] (\x,\y,\z)   -- (\x+1,\y,\z)   -- (\x+1, \y+1, \z)   -- (\x, \y+1, \z) -- cycle;  
      }
      \node at (\x+0.5, \y+0.5, \n) {\n};
 \fi
 }
}
\end{tikzpicture}
\end{minipage}
  \begin{minipage}{0.15\linewidth}  
    \begin{tikzpicture}[x=(220:0.4cm), y=(-40:0.4cm), z=(90:0.4242cm)]
\foreach \m [count=\y] in {{1,1},{1}}{
  \foreach \n [count=\x] in \m {
  \ifnum \n>0
      \foreach \z in {1,...,\n}{
        \draw [fill=blue!30] (\x+1,\y,\z) -- (\x+1,\y+1,\z) -- (\x+1, \y+1, \z-1) -- (\x+1, \y, \z-1) -- cycle;
        \draw [fill=blue!40] (\x,\y+1,\z) -- (\x+1,\y+1,\z) -- (\x+1, \y+1, \z-1) -- (\x, \y+1, \z-1) -- cycle;
        \draw [fill=blue!10] (\x,\y,\z)   -- (\x+1,\y,\z)   -- (\x+1, \y+1, \z)   -- (\x, \y+1, \z) -- cycle;  
      }
      \node at (\x+0.5, \y+0.5, \n) {\n};
 \fi
 }
}
\end{tikzpicture}
\end{minipage}}
\end{center}
\smallskip

 Plane partitions have  also appeared in the study of the Kronecker coefficients in the work of E. Vallejo, see \cite{MR1747064}.
 
We also observe that our results imply that this family of coefficients satisfies the saturation conjecture for the reduced Kronecker coefficients of Kirillov and Klyachko, and that they are weakly increasing. 
Finally we go back to our original aim. In Theorem  \ref{quasipoly}, we show that  the family $\overline{g}_{(k^a),(k^a)}^{(k)}$ is described
by a quasipolynomial of degree $2a-1$  and period dividing the lowest common divisor of $1, 2, \ldots, a+1$.   


\section{Reduced Kronecker coefficients}

In 1938 Murnaghan observed that the Kronecker coefficients $g_{\mu,\nu}^{\lambda}$ always stabilize when 
 we  increase the sizes of the first parts  of the three labelling partitions, see  \cite{MR1507347} and \cite{MR0075213}.    For example, the sequence
$	g_{(k,3,2),(k-1,4,2)}^{(k,2,2,1)} =18, 35, 40, 40,  \ldots$, defined  for $k\ge 5$, has $40$ as its stable value.  

 Denote by $\alpha[n]$ the sequence of integers defined by prepending a first part of size
$n-|\alpha|$ to $\alpha$, then the reduced Kronecker coefficient $ \overline{g}_{\alpha \beta}^\gamma$ is defined to be the stable limit of the sequence $ g_{\alpha[n] \beta[n]}^{\gamma[n]}$. 
Murnaghan noticed that not only each particular sequence $ g_{\alpha[n] \beta[n]}^{\gamma[n]}$ stabilizes,  the Kronecker product $s_{\alpha[n]} \ast s_{\beta[n]}$,  written in the Schur basis, is also stable. Moreover, in \cite[Theorem 1.2]{MR2774644}, it is proved that the Kronecker product $s_{\alpha[n]} \ast s_{\beta[n]}$ stabilizes at $stab(\alpha, \beta) = |\alpha| + |\beta| + \alpha_1 + \beta_1$. 
Let  $stab(\alpha,\beta,\gamma)$ be $\min \{ stab(\alpha,\beta), stab(\alpha, \gamma), stab(\beta, \gamma) \}$.
The symmetry of the Kronecker coefficients implies that if  $n \geq   stab(\alpha,\beta,\gamma)$, then
$
 \overline{g}_{\alpha \beta}^\gamma = g_{\alpha[n] \beta[n]}^{\gamma[n]}.
$

The reduced Kronecker coefficients are interesting objects of their own right.  Littlewood observed that when $|\alpha|+|\beta|=|\gamma|$ they coincide with the Littlewood--Richardson coefficient.  Even if they are believed to be easier to understand than the Kronecker coefficients,  they have be shown to contain  enough information to compute from them the Kronecker coefficients, \cite{MR2774644}.
The Kronecker coefficients  do not satisfy the saturation hypothesis. For example $g_{(n,n), (n,n)}^{(n,n)} $ is equal to $1$ if $n$ is even, and to $0$ otherwise. On the
other hand, both
 Kirillov and Klyachko, in \cite{Klya04} and \cite{MR2105706}, have conjectured that the reduced Kronecker coefficients also satisfy  the saturation hypothesis.
 Recently, the reduced Kronecker coefficients have been used to investigate the rate of grow of the Kronecker coefficients, \cite{BERARM}.


\section{The generating function of a family of reduced Kronecker coefficient}

We compute  the generating function for a special family of reduced Kronecker coefficients. 
 \begin{theorem}\label{ThmRowsKT}
Fix integers $a\geq b\geq 0$. Consider the sequence of reduced Kronecker coefficients $\left\{ \overline{g}^{(k)}_{(k^a), (k^b)} \right\}_{k\geq 0}$.
\begin{itemize}
\item[\textbf{1.}] If $a=b$, the generating function for the reduced Kronecker coefficients $\overline{g}^{(k)}_{(k^a), (k^a)}$ is
\begin{eqnarray*}
\mathcal{F}_{a,a} = \frac{1}{(1-x)\cdot (1-x^2)^2 \cdots  (1-x^a)^2 \cdot (1-x^{a+1})}
\end{eqnarray*}
\item[\textbf{2.}] If $a=b+1$, then $\overline{g}^{(k)}_{(k^a), (k^b)}=1$ for all $k\geq 0$. That is $\mathcal{F}_{b+1,b}=\frac{1}{1-x}$.
\item[\textbf{3.}] If $a > b+1$, $\overline{g}^{(k)}_{(k^a), (k^b)}=0$, except for $k=0$ that it is 1.
\end{itemize}
\end{theorem}
\medskip 

Since $\mathcal{F}_{a,a}$ is the generating function for coloured partitions with parts in $\mathcal{A}= \left\{ \overline{1},2,\overline{2},\dots, a, \overline{a}, \overline{a+1}\right\}$, to prove that Theorem \ref{ThmRowsKT} holds
it suffices to see that $\mathcal{F}_{a,a}$ is also the generating function for the reduced Kronecker coefficients  $\overline{g}^{(k)}_{(k^a), (k^a)}$. We  give an explicit map between  coloured partitions with parts in $\mathcal{A}$, and Kronecker tableaux. 

A SSYT of shape $\lambda / \alpha$ and type $\nu / \alpha$ such that it reverse reading word is an $\alpha-$ lattice permutation is called a \emph{Kronecker tableau} of shape $\lambda / \alpha$ and type $\nu / \alpha$ if either $\alpha_1=\alpha_2$ or $\alpha_1>\alpha_2$ and any one of the following two conditions is satisfied: the number of $1$'s in the second row of $\lambda/\alpha$ is exactly $\alpha_1-\alpha_2$ or the number of $2$'s in the first row of $\lambda/\alpha$ is exactly $\alpha_1-\alpha_2$.

\begin{lemma}[Thm 3.2(a)\cite{MR2264933}]\label{ThmKT}
Let $n$ and $p$ be positive integers such that $n \geq 2p$. 

Let $\lambda=(\lambda_1, \dots, \lambda_{\ell(\lambda)})$ and $\nu$ be partitions of $n$. If $\lambda_1 \geq 2p-1$, then $g_{\lambda, \nu}^{(n-p,p)}= {\displaystyle \sum_{\stackrel{\alpha \vdash p}{\alpha \subseteq
\lambda\cap \nu}}k_{\alpha\nu}^{\lambda}.}$ 
\end{lemma}

The reduced Kronecker coefficient $\overline{g}^{(k)}_{(k^a), (k^a)}$ is equal to the  Kronecker coefficient ${g}^{(k)[N]}_{(k^a)[N], (k^a)[N]}$, when  $N\ge stab(\lambda,\mu,\nu)$.  Lemma \ref{ThmKT},  implies that $\overline{g}^{(k)}_{(k^a), (k^a)}$ is  equal to  the number of Kronecker tableaux with shape $(3k,k^a)/\alpha$ and type $(3k,k^a)/\alpha$, where $\alpha$ is a partition of $k$. 

The following algorithm defines a  bijection between coloured partitions of $k$ with parts in $\mathcal{A}$, and Kronecker tableaux of shape $(3k,k^a)/\alpha$ and type $(3k,k^a)/\alpha$, for $\alpha$ a partition of $k$.

To a  coloured  partition  of $k$, $\beta$, with parts in $\mathcal{A}$,  we associate a Kronecker tableau $T(\beta)$ as follows.  

First, we identify each element of $\mathcal{A}$ to a column of height $a+1$:  if $i \in \{2,3,\dots,a-1,a\}$ then\\

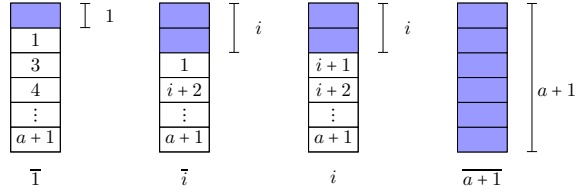
\begin{center}
\scalebox{0.66}{
\begin{tikzpicture}
\draw (0,0.5) rectangle (1,3.5); 
\draw (0,0.5) rectangle (1,1);
\draw (0,1) rectangle (1,1.5);
\draw (0,1.5) rectangle (1,2);
\draw (0,2) rectangle (1,2.5);
\draw (0,2.5) rectangle (1,3);
\draw[fill=blue!40] (0,3) rectangle (1,3.5);
\node at (0.5,0) {$\overline{1}$};
\node at (0.5,0.8) {$a+1$};
\node at (0.5,1.25) {$\vdots$};
\node at (0.5,1.75) {$4$};
\node at (0.5,2.25) {$3$};
\node at (0.5,2.75) {$1$};
\draw[|-|] (1.5,3) -- (1.5,3.5);
\node at (2, 3.25) {$1$};

\draw[fill=blue!40] (9,0.5) rectangle (10,3.5); 
\draw[fill=blue!40] (9,0.5) rectangle (10,1);
\draw[fill=blue!40] (9,1) rectangle (10,1.5);
\draw[fill=blue!40] (9,1.5) rectangle (10,2);
\draw[fill=blue!40] (9,2) rectangle (10,2.5);
\draw[fill=blue!40] (9,2.5) rectangle (10,3);
\draw[fill=blue!40] (9,3) rectangle (10,3.5);
\node at (9.5,0) {$\overline{a+1}$};
\draw[|-|] (10.5,0.5) -- (10.5,3.5);
\node at (11, 1.75) {$a+1$};

\draw (3,0.5) rectangle (4,3.5); 
\draw (3,0.5) rectangle (4,1);
\draw (3,1) rectangle (4,1.5);
\draw (3,1.5) rectangle (4,2);
\draw (3,2) rectangle (4,2.5);
\draw[fill=blue!40] (3,2.5) rectangle (4,3);
\draw[fill=blue!40] (3,3) rectangle (4,3.5);
\node at (3.5,0) {$\overline{i}$};
\node at (3.5,0.8) {$a+1$};
\node at (3.5,1.25) {$\vdots$};
\node at (3.5,1.75) {$i+2$};
\node at (3.5,2.25) {$1$};
\draw[|-|] (4.5,2.5) -- (4.5,3.5);
\node at (5, 3) {$i$};

\draw (6,0.5) rectangle (7,3.5);  
\draw (6,0.5) rectangle  (7,1);
\draw (6,1) rectangle (7,1.5);
\draw (6,1.5) rectangle (7,2);
\draw (6,2) rectangle (7,2.5);
\draw[fill=blue!40] (6,2.5) rectangle (7,3);
\draw[fill=blue!40] (6,3) rectangle (7,3.5);
\node at (6.5,0) {$i$};
\node at (6.5,0.8) {$a+1$};
\node at (6.5,1.25) {$\vdots$};
\node at (6.5,1.75) {$i+2$};
\node at (6.5,2.25) {$i+1$};
\draw[|-|] (7.5,2.5) -- (7.5,3.5);
\node at (8, 3) {$i$};
\end{tikzpicture}}
\end{center}

Note that it is always possible to order the  columns  corresponding to the parts of $\beta$ in such a way that we obtain a semistandard Young tableau. If we write $\beta$ as $\left( \overline{1}^{m_{\overline{1}}} \overline{2}^{m_{\overline{2}}} 2^{m_2} \dots \overline{a+1}^{m_{\overline{a+1}}} \right)$, then $m_i$ will denote the number of times that the column $i$ appears in the tableau that we are building. 

We read the partition  $\alpha$ from our semistandard Young tableau by counting the number of blue boxes in each row: $\alpha_{a+1}= m_{\overline{a+1}}$, $\alpha_i = \alpha_{i+1} +m_i + m_{\overline{i}}$ for $i=2,\dots, a$, and $\alpha_1= \alpha_2 +m_{\overline{1}}$.

This semistandard Young tableau will be the first columns on the left side of $T(\beta)$. We build the rest of the Kronecker tableau of shape $(3k,k^a)/\alpha$ as follows: complete the $i-$th row with $i$'s, for $i=2,\dots,a+1$, and complete the first row with the remaining numbers of the type $(3k,k^a) / \alpha$ in weakly increasing order from left to right. 
For instance, \\

\begin{center}
   \begin{minipage}{0.4\linewidth}  
  Kronecker tableau corresponding to \\
   $\lambda= \nu=(9,3,3,3)$ and $\alpha=(2,1)$.\\
  It is obtained by our algorithm taking \\
  $a=3$ and $\beta=(\overline{2},\overline{1})$. \\
\end{minipage}
{\small
\begin{minipage}{0.5\linewidth}
\begin{center}
\begin{tikzpicture}
\draw (0,0.5) rectangle (1.5,2.5);
\draw (0,1) -- (1.5,1);
\draw (0,1.5) -- (1.5,1.5);
\draw (0,2) -- (4.5,2);
\draw (0,2.5) -- (4.5,2.5);
\draw (0.5,0.5) -- (0.5,2.5);
\draw (1,0.5) -- (1,2.5);
\draw (1.5,0.5) -- (1.5,2.5);
\draw (2,2)-- (2,2.5);
\draw (2.5,2)-- (2.5,2.5);
\draw (3,2)-- (3,2.5);
\draw (3.5,2)-- (3.5,2.5);
\draw (4,2)-- (4,2.5);
\draw (4.5,2)-- (4.5,2.5);

\draw[fill=blue!40] (0,1.5) rectangle (0.5,2);
\draw[fill=blue!40] (0,2) rectangle (0.5,2.5);
\draw[fill=blue!40] (0.5,2) rectangle (1,2.5);

\node at (0.25,0.75) {4};
\node at (0.25,1.25) {1};
\node at (0.75,0.75) {4};
\node at (0.75,1.25) {3};
\node at (0.75,1.75) {1};
\node at (1.25,0.75) {4};
\node at (1.25,1.25) {3};
\node at (1.25,1.75) {2};
\node at (1.25,2.25) {1};
\node at (1.75,2.25) {1};
\node at (2.25,2.25) {1};
\node at (2.75,2.25) {1};
\node at (3.25,2.25) {1};
\node at (3.75,2.25) {2};
\node at (4.25,2.25) {3};
\end{tikzpicture}
\end{center}
  \end{minipage}
  }
\end{center}

\vspace{0.3cm} 
This map  is well-defined and bijective. For the  other two cases we  show that there is either only one  Kronecker tableau that satisfy all requirements, or none at all. 


\section{Plane Partitions and Reduced Kronecker Coefficients}

In this section we establish a link between our family of reduced Kronecker coefficients and   plane partitions. 
 A \emph{plane partition} is a finite subset $\mathcal{P}$ of positive integer lattice points, $\{(i,j,k)\} \subset \mathbb{N}^3$, such that if $(r,s,t)$ lies in $\mathcal{P}$ and if $(i,j,k)$ satisfies $1\leq i \leq r$, $1\leq j \leq s$ and $1\leq k\leq t$,  then $(i,j,k)$ also lies in $\mathcal{P}$. 
 Let $\mathcal{B}(r,s,t)$ be the set of plane partitions fitting in a $r\times s$ rectangle and biggest part equals to $t$. 
That is, 
$$  \mathcal{B} (r,s,t)= \left\{ (i,j,k)\ |\ 1\leq i \leq r,\ 1\leq j \leq s,\ 1\leq k\leq t \right\}.$$

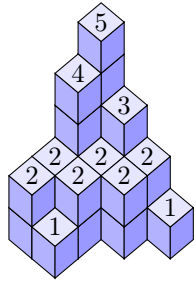
\begin{center}
  \begin{minipage}{0.3\linewidth}  
  {\small
    \begin{tikzpicture}[x=(220:0.4cm), y=(-40:0.4cm), z=(90:0.4242cm)]
\foreach \m [count=\y] in {{5,4,2,2},{3,2,2,1},{2,2},{1}}{
  \foreach \n [count=\x] in \m {
  \ifnum \n>0
      \foreach \z in {1,...,\n}{
        \draw [fill=blue!30] (\x+1,\y,\z) -- (\x+1,\y+1,\z) -- (\x+1, \y+1, \z-1) -- (\x+1, \y, \z-1) -- cycle;
        \draw [fill=blue!40] (\x,\y+1,\z) -- (\x+1,\y+1,\z) -- (\x+1, \y+1, \z-1) -- (\x, \y+1, \z-1) -- cycle;
        \draw [fill=blue!10] (\x,\y,\z)   -- (\x+1,\y,\z)   -- (\x+1, \y+1, \z)   -- (\x, \y+1, \z) -- cycle;  
      }
      \node at (\x+0.5, \y+0.5, \n) {\n};
 \fi
 }
}
\end{tikzpicture}}
\end{minipage}
\begin{minipage}{0.25\linewidth}
A plane partition  in $\mathcal{B}(4,4,5)$.
Its two-dimensional array is\\
\begin{center}
\begin{tikzpicture}
\draw (0,0) rectangle (2,2);
\draw (0,0) rectangle (0.5,0.5);
\node at (0.25,0.25) {2};
\draw (0.5,0) rectangle (1,0.5);
\node at (0.75,0.25) {1};
\draw (0,0.5) rectangle (0.5,1);
\node at (0.25,0.75) {2};
\draw (0.5,0.5) rectangle (1,1);
\node at (0.75,0.75) {2};
\draw (0,1) rectangle (0.5,1.5);
\node at (0.25,1.25) {4};
\draw (0.5,1) rectangle (1,1.5);
\node at (0.75,1.25) {2};
\draw (1,1) rectangle (1.5,1.5);
\node at (1.25,1.25) {2};
\draw[fill=blue!20] (0,1.5) rectangle (0.5,2);
\node at (0.25,1.75) {5};
\draw (0.5,1.5) rectangle (1,2);
\node at (0.75,1.75) {3};
\draw (1,1.5) rectangle (1.5,2);
\node at (1.25,1.75) {2};
\draw (1.5,1.5) rectangle (2,2);
\node at (1.75,1.75) {1};
\end{tikzpicture}
\end{center}
  \end{minipage}
\end{center}

\vspace{0.3cm} We show that the generating function for the reduced Kronecker coefficients, obtained in Theorem \ref{ThmRowsKT}, coincides with the  classical generating functions for  plane partitions.
\begin{theorem}[P. MacMahon, \cite{MR2417935}]\label{ThmPPIni}
 The generating function for plane partitions that are subsets of $\mathcal{B}(r,s,t)$ is given by
$$
  pp_t(x;r,s)= \prod_{i=1}^r \prod_{j=1}^s \frac{1-x^{i+j+t-1}}{1-x^{i+j-1}}
$$
\end{theorem}
Since $\frac{1-x^{l+k}}{1-x^l}= 1+x^l+\cdots + x^{\left\lfloor \frac{k}{l}\right\rfloor l} + \mathcal{O}(k+1)$,  for all $l\geq 1$, this generating function for plane partitions  can be rewritten  as a generating function computed resembling the one appearing in  Theorem \ref{ThmRowsKT}.
 \begin{lemma}\label{ThmPP}
Let $r=\min(a,l)$ and $s=\max(a,l)$. Then, the generating function for the plane partitions fitting inside an $l\times a$ rectangle is
$$
 \prod_{j=r}^{s} \left( \frac{1}{1-x^j} \right)^{r} \cdot \prod_{i=1}^{r-1}\left(\frac{1}{1-x^i}\right)^i \left( \frac{1}{1-x^{s + i}}\right)^{r-i}
$$
  \end{lemma}
From this we easily conclude that the family of reduced Kronecker coefficients $\overline{g}_{(k^a), (k^a)}^{(k)}$ counts plane partitions inside a certain rectangle. That is, 
 \begin{theorem}\label{ThmPPRows}
 The reduced Kronecker coefficient $\overline{g}_{(k^a), (k^a)}^{(k)}$ counts the number of plane partitions of $k$ fitting inside a $2\times a$ rectangle.  
 \end{theorem}


 \section{Consequences} 
 
 {\bf 4.1 Saturation Hypothesis:} 
Let denote by $\{ C(\alpha^1,\dots, \alpha^n) \}$ any family of coefficients depending on the partitions $\alpha^1,\dots, \alpha^n$. The family  $\{ C(\alpha^1,\dots, \alpha^n) \}$ satisfies the \emph{saturation hypothesis} if the conditions 
$ C(\alpha^1,\dots, \alpha^n) > 0 $ and $ C(s\cdot \alpha^1, \dots, s\cdot \alpha^n) > 0 $ for all $s>1$ are equivalent, where $s\cdot \alpha = (s\alpha_1, s\alpha_2,\dots)$.  The  Littlewood-Richardson coefficients  satisfy the saturation hypothesis, as shown by Knutson and Tao in \cite{MR1671451}. On the other hand, the Kronecker coefficients are known not to satisfy it.

In \cite{Klya04} and \cite{MR2105706}, Kirillov and Klyachko  have conjectured that the  reduced Kronecker coefficients satisfy the saturation hypothesis.
 From the combinatorial interpretation for the reduced Kronecker coefficients  $\overline{g}_{(k^a), (k^a)}^{(k)}$ in terms of plane partitions we verify their conjecture for our family of coefficients.
 \begin{corollary}
The saturation hypothesis holds for the coefficients $\overline{g}_{(k^a), (k^a)}^{(k)}$. In fact, $\overline{g}_{((sk)^a), ((sk)^a)}^{(sk)} >0$ for all $s\geq 1$. Moreover, the sequences of coefficients obtained by, either fixing $k$ or $a$, and then  letting the other parameter grow are weakly increasing. 
\end{corollary}

\medskip

{\bf 4.2 Quasipolynomiality:} In Theorem \ref{ThmRowsKT} we computed the generating function  $\mathcal{F}_{a,b}$ for the reduced Kronecker coefficients. In this section we study the implications of this calculation. We concentrate on the non-trivial case, $a=b$. 
 
 \begin{theorem}\label{quasipoly}
  Let $\mathcal{F}_{a}=  \mathcal{F}_{a,a}$ be the generating function for  the reduced Kronecker coefficients $\overline{g}^{(k)}_{(k^a), (k^a)}$.

Let $\ell$ be the lowest common multiple of $1, 2, \ldots, a, a+1$.
\begin{enumerate}

\item The generating function $\mathcal{F}_{a}$ can be rewritten as
\[
\mathcal{F}_{a}= \frac{P_a(x)}{(1-x^{\ell})^{2a}}
\]
where $P_a(x)$ is a product of cyclotomic polynomials, and $\deg(P_a(x)) =2\ell a - (a+2)a< 2a\ell-1$. 

\item The polynomial $P_a$ is the generating function for coloured partitions with parts in $\{1, 2, \bar 2, 3, \bar 3, \ldots, a, \bar a, a+1 \}$, where parts   $j$ and $\bar j$ appear with multiplicity  less than $\ell/j$ times. 

\item The coefficients of $P_a$ are positive and  palindrome, but in general are not a concave sequence.

\item  The  coefficients $\overline{g}^{(k)}_{(k^a), (k^a)}$ are described by a quasipolynomial of degree $2a-1$ and period dividing $\ell$. In fact, we have checked that the period is exactly $l$ for $a$ less than 10.

\item The  coefficients $\overline{g}^{(k)}_{(k^a), (k^a)}$  satisfy a formal  reciprocity  law: $x^{a(a+2)}\mathcal{F}_a(x)= \mathcal{F}_a(\frac{1}{x})$. 

\end{enumerate}
 \end{theorem}
This theorem is shown using Proposition 4.13 of  \cite{BecSan}.\\

\begin{example}
Let  $\Phi_i$ be the $i^{th}$ cyclotomic polynomials. From the well-known identity $(x^n-1)= \prod_{i | n} \Phi_i$, we express  $P_a$ as a product of cyclotomic polynomials. For example, $P_{2}=\Phi_2^2\Phi_3^3\Phi_6^4$ and $P_3=\Phi_2^3\Phi_3^4\Phi_4^5\Phi_6^6\Phi_{12}^6$. It can be easily seen that $\Phi_1$ never appears in these expansions, so the polynomials are palindromes.
 \end{example} 
 \smallskip
 
 \begin{example}
The coefficients $\overline{g}^{(k)}_{(k^2), (k^2)}$ are given by the quasipolynomial of degree $3$ and period $6$. 
\begin{align*}
\overline{g}^{(k)}_{(k^2), (k^2)} =
 \left\{
{\small{
 \begin{array}{ll}
 1/72 k^3 +1/6 k^2+\phantom{1}2/3  k +\phantom{5}1 \phantom{1}                            & \text{ if } k \equiv 0 \mod 6\\
 1/72 k^3 +1/6 k^2+13/24 k +5/18             & \text{ if } k \equiv 1 \mod 6\\
 1/72 k^3 +1/6 k^2+\phantom{1}2/3  k +8/9\phantom{1}              & \text{ if } k \equiv 2 \mod 6\\
 1/72 k^3 +1/6 k^2+ 13/24  k +1/2\phantom{1}                       &\text{ if } k \equiv 3 \mod 6\\
 1/72 k^3 +1/6 k^2+\phantom{1}2/3  k +7/9\phantom{1}                             &\text{ if } k \equiv 4 \mod 6\\
 1/72 k^3 +1/6 k^2+13/24  k +7/18                        &\text{ if } k \equiv 5 \mod 6
 \end{array}
  }}
 \right.
 \end{align*}
 \end{example} 
 
 \medskip

{\bf 4.3 On the grow of the Kronecker coefficients:} 
Murnaghan famously  observed that the sequences obtained by adding cells to the first parts of the partitions indexing a Kronecker coefficients are eventually constant.
Fix three arbitrary partitions. In \cite{BERARM}, it is shown that the sequences  obtained by adding cells to the first parts of the three partitions indexing a reduced Kronecker coefficients are described by a linear quasipolynomials of period 2. 

These sequences can be interpreted as adding cells to the second rows of the partitions indexing a Kronecker coefficient, while keeping their fist parts very long in comparison.
This result is an extension of Murnaghan's  observation to the other rows of the partition.

  An interesting question is then to describe what happens when we add cells to arbitrary rows of the partitions indexing a Kronecker (and reduced Kronecker) coefficient. The results presented in this note can be seen as a contribution to this investigation. We have shown that for any $a$, the sequence $\bar g^{k}_{{k^a},{k^a}}$  is  described by a quasipolynomial of degree $2a-1$ and period dividing $\ell$.

 On the other hand, in  \cite{BERARM} it is also shown that when we fix three partitions, and start adding cells to their first columns, the sequences obtained are eventually constant.
 From our combinatorial interpretation for $\bar g^{k}_{{k^a},{k^a}}$ in terms of plane partitions of $k$ fitting a $2 \times a$ rectangle we obtain that, if we fix $k$ and let $a$ grow, the sequences obtained are always going to be eventually constant. Moreover,  for any $a>k$ the $\bar g^{k}_{{k^a},{k^a}}$ counts the number of such plane partitions of $k$ with at most two parts.


\bibliographystyle{alpha}
\bibliography{Referencias_Articulos}
\label{sec:biblio}



\end{document}